\documentclass[reqno,11pt]{amsart}

\usepackage[all]{xy}

\numberwithin{equation}{section}

\theoremstyle{plain}
\newtheorem{lemma}{Lemma}[section]
\newtheorem{theor}[lemma]{Theorem}
\newtheorem{prop}[lemma]{Proposition}
\newtheorem{cor}[lemma]{Corollary}

\newtheorem{claim}[lemma]{Claim}

\theoremstyle{definition}

\newtheorem{exa}[lemma]{Example}

\theoremstyle{remark}
\newtheorem{rem}[lemma]{Remark}

\newcommand{\pic}{{\rm Pic\thinspace}} 
\newcommand{\Pic}{{\rm Pic\thinspace}}
\newcommand{\bs}{{\rm Bs}}
\newcommand{\Bs}{{\rm Bs}}
\newcommand{\Cr}{{\rm Cr}}

\newcommand{\p}{\mathbb{P}}

\newcommand{\z}{\mathbb{Z}}
\newcommand{\f}{\mathbb{F}}

\newcommand{\co}{{\mathcal O}}
\newcommand{\cl}{{\mathcal L}}

\newcommand{\cm}{{\mathcal M}}
\newcommand{\cn}{{\mathcal N}}

\newcommand{\cs}{{\mathcal S}}

\newcommand{\ci}{{\mathcal I}}

\newcommand{\et}{\tilde{e}}
\newcommand{\hti}{\tilde{h}}

\newcommand{\vs}{\vspace{2mm}}
\newcommand{\svs}{\vspace{1mm}}
\newcommand{\nvs}{\vspace{-2mm}}
\newcommand{\bnvs}{\vspace{-4mm}}

\newcommand{\map}{\rightarrow}
\newcommand{\rmap}{\dashrightarrow}

\newcommand{\vrk}{\hfill $\Box$}

\begin{document}

\title[Base locus of linear systems on blowings-up of $\p^3$]{Base locus of linear systems on the blowing-up of $\p^3$ along at most 8 general points}
\author{Cindy De Volder}
\address{
Department of Pure Mathematics and Computeralgebra,
Galglaan 2, \newline B-9000 Ghent, Belgium}
\email{cdv@cage.ugent.be}
\thanks{The first author is a Postdoctoral Fellow of
the Fund for Scientific Research-Flanders (Belgium)
(F.W.O.-Vlaanderen)}

\author{Antonio Laface}
\address{
Dipartimento di Matematica, Universit\`a degli Studi di Milano,
Via Saldini 50, \newline 20100 Milano, Italy }
\email{antonio.laface@unimi.it}
\thanks{The second author would like to acknowledge the support of the MIUR of the
Italian Government in the framework of the National Research
Project ``Geometry in Algebraic Varieties'' (Cofin 2002)}

\keywords{Base locus, linear systems, fat points, projective space}
\subjclass{14C20}
\begin{abstract}
Consider a (non-empty) linear system of surfaces of degree $d$ in $\p^3$ 
through at most 8 multiple points in general position
and let $\cl$ denote the corresponding complete linear system
on the blowing-up $X$ of $\p^3$ along those general points.
Then we determine the base locus of such linear systems $\cl$ on $X$.
\end{abstract}
\maketitle

\section{Introduction}

In this paper we wok over an algebraically closed field of characteristic $0$.

Let $P_1, \ldots, P_r$ be general points of the $n$-dimensional projective space $\p^n$
and
choose some integers $m_1 , \ldots, m_r$.
Consider the linear system $\cl'$ of hypersurfaces of degree $d$ in $\p^n$ having 
multiplicities at least $m_i$ at $P_i$, for all $i=1,\ldots, r$.
Let $X$ denote the blowing-up of $\p^n$ alongs $P_1, \ldots, P_r$,
and let $\cl$ denote the complete linear system on $X$ corresponding to $\cl'$.

A point $Q \in X$ is called a base point of $\cl$ if 
$Q \in D$ for every divisor $D \in \cl$.
The scheme-theoretical union of all base points of $\cl$
is called the base locus of $\cl$.

In case $n=2$ (i.e. if $X$ is a rational surface obtained 
by blowing-up $\p^2$ along $r$ general points)
the dimension, base locus and other properties of
linear systems $\cl$ has been widely studied
(see e.g.~\cite{CD}, \cite{CM1}, \cite{CM2}, \cite{AG}, \cite{AH}).

In case $n=3$, i.e. if $X$ is a rational threefold obtained by
blowing-up $\p^2$ along $r$ general points,
very little is known.
If $r \leq 8$, the dimension of $\cl$ can be determined using the results from~\cite{DL}.
In this paper (for $n=3$ and $r \leq 8$), we will completely describe the base locus of 
$\cl$ on $X$.

In sections~\ref{prel} to~\ref{hirz} we state some preliminaries and notation.
The main results are formulated in section~\ref{results}
and the last three sections contain their proofs.

\section{Preliminaries}\label{prel}

Let $P_1, \ldots, P_8$ be general points on $\p^3$,
let $X$ denote the blowing-up of $\p^3$ along these 8 points,
denote the projection map by $\pi: X \map \p^3$ and
let $E_i$ be the exceptional divisor corresponding to $P_i$.

By $\cl_3(d ; m_1,\ldots, m_r)$, with $r \leq 8$, we denote the complete linear system on $X$
corresponding to the invertible sheaf
$\pi^*(\co_{\p^3}(d)) \otimes \co_X(-m_1 E_1 - \cdots - m_r E_r)$;
i.e. the complete linear system corresponding to
the linear system of hypersurfaces of degree $d$
with multiplicities at least $m_i$ at $P_i$.
Similarly, by $\cl_3(d ; m_1^{r_1},\ldots, m_s^{r_s})$
(with $r_1 + \cdots + r_s \leq 8$), we denote the complete linear system on $X$
corresponding to the linear system of hypersurfaces of degree $d$
with $r_j$ points of multiplicities at least $m_j$. 

With $\langle h,e_1,\ldots,e_r\rangle$ we
denote a basis of ${\rm\bf A}^2(X)$ where $h$ is the pull-back of
a class of a general line in $\p^3$ and $e_i$ is the class of a
line on $E_i$.
The notation
$\ell=\ell_3(\delta,\mu_1,\ldots,\mu_r)$ indicates the set of 
the strict transforms of all curves in
$\p^3$ of degree $\delta$ through $r$ points of multiplicity
$\mu_1,\ldots,\mu_r$ or equivalently all curves in
$| \delta h-\sum_{i=1}^r\mu_i e_i |$ on $X$.

For $1 \leq i < j \leq 8$, we denote the strict transform of the line
through $P_i$ and $P_j$ by $\ell_{i,j}$.

We say a class $\cl_3(d ; m_1,\ldots, m_r)$ is in standard form if
$m_1 \geq \cdots \geq m_r \geq 0$ and $2d \geq m_1 + m_2 + m_3 + m_4$.
In \cite[Proposition 2.2]{DL} we prove the following

\begin{lemma}\label{standard}
A linear system $\cl = \cl_3(d ; m_1,\ldots, m_r)$ is in standard form if and only if
$\cl = \cs + \sum_{i=4}^a c_i \cs_i$
with $c_i \in \z_{\geq 0}$, $\cs_i = \cl_3(2; 1^i)$
and $\cs=\cl_3(d-2m_4,m_1-m_4,m_2-m_4,m_3-m_4)$. \vrk
\end{lemma}

For all $1 \leq i \leq 8$, let $Q_i$ be a general element of $\cs_i (= \cl_3(2; 1^i))$.
Then $Q_i$ is the blowing-up of $\bar{Q}_i$, a general quadric hypersurface  
in $\p^3$ through
the points $P_1, \ldots, P_i$, 
along those $i$ points.
Also $\Pic Q_i = \langle f_1, f_2, e_1, \ldots, e_i \rangle$,
with $f_1$ and $f_2$ the pullbacks of the two rulings on $\bar{Q}_i$
and $e_1, \ldots, e_i$ the exceptional curves.
By $\cl_{Q_i}(a, b ; m_1 , \ldots ,m_i)$ we denote 
the complete linear system $|a f_1 + b f_2 - m_1 e_1 - \ldots - m_i e_i |$,
and, as before, if some of the multiplicities are the same, we also use
the notation $\cl_{Q_i}(a, b ; m_1^{n_1} , \ldots ,m_r^{n_r})$.

Let $B_j$ be the blowing-up of $\p^2$ along $j$ general points,
then $\Pic B_j = \langle h , e'_1 , \ldots, e'_j \rangle$,
with $h$ the pullback of a line and $e'_l$ the exceptional curves.
By $\cl_2( d ; m_1 , \ldots ,m_j)$ we denote 
the complete linear system $|d h - m_1 e'_1 - \ldots - m_j e'_j |$.
And again, as before, if some of the multiplicities are the same, we also use
the notation $\cl_{2}(d ; m_1^{n_1} , \ldots ,m_r^{n_r})$.

On $B_j$, a system $\cl_2( d ; m_1 , \ldots ,m_j)$ is said to be in standard form if 
$d \geq m_1 + m_2 + m_3$ and $m_1 \geq m_2 \geq \cdots \geq m_j \geq 0$;
and it is called standard it there exists a base 
$\langle \hti , \et_1 , \ldots, \et_j \rangle$ 
of $\Pic B_j$
such that 
$\cl_2( d ; m_1 , \ldots ,m_j) = |\tilde{d} \hti - \tilde{m}_1 \et_1 - \ldots - \tilde{m}_j \et_j |$
is in standard form.

As explained in~\cite[\S 6]{DL}, the blowing-up $Q_i$ of the quadric along $i$ general points
can also be seen as a blowing-up of the projective plane
along $i +1$ general points,
and 
$$\cl_{Q_i}(a, b ; m_1 , \ldots ,m_i) = \cl_2( a+b-m_1 ; a-m_1 , b-m_1 , m_2 ,\ldots, m_i).$$
In particular the anticanonical class $-K_{Q_i}$ contains an irreducible
reduced divisor which we denote by $D_{Q_i}$.

\section{Cubic Cremona transformation}

The cubic Cremona transformation on $\p^3$, 
whose associated rational map is given by  
\begin{equation}\label{cubic}
\begin{array}{rccl}
\Cr: &  \p^3  &  \rmap  &  \p^3  \\
     & (x_0:x_1:x_2:x_3) & \mapsto  & (x_0^{-1}:x_1^{-1}:x_2^{-1}:x_3^{-1}),
\end{array}
\end{equation}    
induces an action on $\Pic X$, resp. on ${\rm\bf A}^2(X)$,
as stated in the following two propositions
(see~\cite{LU} for a proof of both results). 

\begin{prop}\label{cre surf}
Assuming the Cremona transformation~(\ref{cubic}) uses the points $P_1,\ldots, P_4$,
its induced action on
$\cl = \allowbreak \cl_3(d,\allowbreak m_1,\allowbreak \ldots,m_r)$
is given by
\begin{eqnarray}\label{form cre surf}
\Cr(\cl) & := & \cl_3(d+k,m_1+k,\ldots, m_4+k,m_5,\ldots,m_r),
\end{eqnarray}
where $k=2d-\sum_{i=1}^4m_i$.    \vrk
\end{prop}   

\begin{prop}\label{cre curve}
Assuming the Cremona transformation~(\ref{cubic}) uses the points $P_1,\ldots, P_4$,
its induced action on
$\ell = \allowbreak \ell_3(\delta,\allowbreak \mu_1,\allowbreak \ldots,\mu_r)$,
with $\ell$ skew to the $\ell_{i,j}$ for $1 \leq i < j \leq 4$, 
is given by
\begin{eqnarray}\label{form cre curve}
\Cr(\ell) & := & \ell_3(\delta+2h,\mu_1+h,\ldots,       
\mu_4+h,\mu_5,\ldots,\mu_r),
\end{eqnarray}
where $h=\delta-\sum_{i=1}^4\mu_i$. 
Moreover, under the same assumption,
for all 
$1 \leq i < j \leq 4$,
we have that
$\Cr(\ell_{i,j}) = \ell_{u,v}$, with $\{i,j,u,v\} = \{1,2,3,4\}$.
\vrk
\end{prop}

\begin{rem}
It follows immediately from the previous propositions that the
Cremona transformation fixes $\cs_i$ (for $4 \leq i \leq 8$)
and $K_{Q_8} (= \ell_3(4 ; 1^8))$,
i.e. $\Cr (\cs_i) = \cs_i$ and $\Cr (K_{Q_8}) = K_{Q_8}$.
Moreover $\Cr (\cl) . K_{Q_8} = \cl . K_{Q_8}$.
\end{rem}

\begin{rem}\label{cr on X}
It can be proved (see~\cite{LU}) that the cubic Cremona transformation on $X$,
is obtained by blowing-up the strict transforms of the six edges of the tetrahedron
through the four points used by the cubic Cremona transformation,
and blowing down along the other rulings of the exceptional quadrics.
This implies in particular that the cubic Cremona transformation 
is not just a base change of $\Pic X$.
\end{rem}

Let $Y$ denote the blowing-up of $X$ along 
the $l_{1,2}, l_{1,3}, l_{1,4}, l_{2,3}, l_{2,4}$ and $l_{3,4}$.
Then 
$$ \pic Y = \langle H, E_1, \ldots, E_8, E_{1,2} , \ldots, E_{3,4} \rangle $$
where $H$ is the pull-back of a plane in $\p^3$,
$E_i$ is the pull-back of $E_i$ on $X$ (for all $1 \leq i \leq 8$) and
$E_{i,j}$ is the exceptional quadric corresponding to $l_{i,j}$ 
(for all $1 \leq i < j \leq 4$).

On $Y$ the Cremona transformation using the points $P_1, \ldots, P_4$, 
is then nothing else than a base change for $\pic Y$.
In particular, in~\cite{LU}, it is shown that 
\begin{equation}\label{pic Y}
\begin{aligned}
\pic Y & =  \langle H, E_1, \ldots, E_8, E_{1,2} , \ldots, E_{3,4} \rangle \\
       & =  \langle H', F_1, \ldots, F_4, E_5, \ldots, E_8, F_{1,2} , \ldots, F_{3,4} \rangle
\end{aligned}
\end{equation}
with
\begin{equation}\label{F i,j}
\begin{gathered}
\!\!\!\!\!\!\!\!\!\!\!\!\!\!\!\!\!\!\!\!\!\!\!\!
	H'  = \Cr(H) = 3 H \, - \sum_{i=1}^4 2 E_i \, - \!\! \sum_{1 \leq i < j \leq 4} E_{i,j} \!\! \\
\!\!\!\!\!\!\!\!\!\!\!\!\!\!\!\!\!\!\!\!\!\!\!\!
	F_k  = \Cr (E_k) = H \, - \!\!\! \sum_{\substack{1 \leq j \leq 4 \\ j \neq k}} E_j \, - \!\!\!
		\sum_{\substack{1 \leq i<j \leq 4 \\ i,j \neq k}} E_{i,j} \\
\!\!\!\!\!\!\!\!\!\!\!\!\!\!\!\!\!\!\!\!\!\!\!\!
	F_{i,j}  = \Cr (E_{i,j}) = E_{k,l} \mbox{ with } \{i,j,k,l\} = \{1,2,3,4\}.
			 \!\!\!\!\!\!\!\!\!\!\!\!\!\!\!\!\!\!\!\!\!\!  
\end{gathered}
\end{equation}

It follows immediately from these formula that
\begin{multline}\label{cr L on Y}
|dH \, - \!\!\! \sum_{1 \leq i \leq 4} m_i E_i \,\, - 
		\!\!\!\! \sum_{1 \leq i < j \leq 4} \!\!\!\!  m_{i,j} E_{i,j}| \\
 	=  |(d+s)H' \, - \!\!\!\ \sum_{1 \leq i \leq 4} (m_i +s) F_i \,\, 
		- \!\!\!\!\!\!\!\!\!\!\!\!\!\!
		\sum_{\substack{1 \leq i<j \leq 4 \\ \{i,j,k,l\} = \{1,2,3,4\}} } \!\!\!\!\!\!\!\!\!\!\!\!
			(d - m_k - m_l + m_{k,l}) F_{i,j}|
\end{multline}

Similarly, for ${\rm\bf A}^2(Y)$, we have (see~\cite{LU})
\begin{equation}\label{A^2(Y)}
\begin{aligned}
{\rm\bf A}^2(Y)	& =  \langle h, e_1, \ldots, e_8, e_{1,2} , \ldots, e_{3,4} \rangle \\
       		& =  \langle h', f_1, \ldots, f_4, e_5, \ldots, e_8, f_{1,2} , \ldots, f_{3,4} \rangle
\end{aligned}
\end{equation}
with $h$ the pull-back of a line in $\p^3$, $e_i$ the class of a line in $E_i$,
$e_{i,j}$ the vertical ruling of $E_{i,j}$ and
\begin{equation}\label{f i,j}
\begin{aligned}
	h' 	& = \Cr(h) = 3 h \, - \sum_{i=1}^4 e_i \\\
	f_k 	& = \Cr (e_k) = 2 h \, - \!\!\! \sum_{\substack{1 \leq j \leq 4 \\ j \neq k}} e_j \\
	f_{i,j} & = \Cr (e_{i,j}) = h + e_{k,l} - e_k - e_l \mbox{ with } \{i,j,k,l\} = \{1,2,3,4\}.
\end{aligned}
\end{equation}

For the rest of this paper, 
we will use the sheaf notation (e.g. $\pi^*(\co_{\p^3}(d)) \otimes \co_X(-m_1 E_1)$) 
as well as the linear system notation (e.g. $|dH - m_1 E_1|$)
for both purposes, it should be clear from the context which one is intended.

\section{(-1)-curves on $X$}

A curve $C \in \ell = \ell_3(\delta,\allowbreak \mu_1,\allowbreak \ldots,\mu_r)$ 
is called a $(-1)$-curve if $\ell$ is obtained by
applying a finite set of cubic Cremona transformations on the system
$\ell_3(1,1^2)$.

For all $a \in \z_{\geq 0}$ and $b\stackrel{\neq}{,} c \in \{1, \ldots,8\}$,
let 
$$
\delta_{i;b,c} = 
    \begin{cases}
        0  &  \mbox{if } i \notin \{b,c\} \\
        1  &  \mbox{if } i \in \{b,c\};
    \end{cases}
$$
and
$$
{\mathcal C}_a^{b,c} = 
    \begin{cases}
        \ell_3(2a+1 ; \frac{a}{2} + \delta_{1;b,c} , \frac{a}{2} + \delta_{2;b,c} , \cdots ,
                \frac{a}{2} + \delta_{8;b,c})  &   \mbox{if } a \mbox{ is even} \\
        \ell_3(2a+1 ; \frac{a+1}{2} - \delta_{1;b,c} , \frac{a+1}{2} - \delta_{2;b,c} , \cdots , 
                 \frac{a+1}{2} - \delta_{8;b,c})  &   \mbox{if } a \mbox{ is odd}.               
    \end{cases}
$$                                                                                    

\begin{lemma}\label{-1 = C_a^{b,c}}
A curve $C \in \ell$ on $X$ is a (-1)-curve if and only if 
there exists $a \in \z_{\geq 0}$ and $b\stackrel{\neq}{,} c \in \{1, \ldots,8\}$
such that $\ell = {\mathcal C}_a^{b,c}$.\nvs
\end{lemma}

\begin{proof}
First of all, note that $\ell_{i,j} = {\mathcal C}_0^{i,j}$.
So all ${\mathcal C}_0^{i,j}$ are classes of $(-1)$-curves.
To simplify notation, we now assume that $i=1$ and $j=2$,
and by $B$ we denote the set of the four indices of the points used for
the transformation~(\ref{cubic}).
To determine $\Cr(\ell_{1,2})$
we distinguish three cases
\begin{itemize}
	\item[(a)] $P_1$ and $P_2$ $\in B$.
		Without loss of generality we may assume 
		that the transformation~(\ref{cubic}) uses the points $P_1,\ldots,P_4$,
		i.e. that $B = \{1,2,3,4\}$.
		So, by proposition~\ref{cre curve} we obtain 
		that $\Cr(\ell_{1,2}) = \ell_{3,4}$, i.e.
		$\Cr({\mathcal C}_0^{1,2}) = {\mathcal C}_0^{3,4}$.
	\item[(b)] $P_2 \in B$ and $P_1 \notin B$.
		Then we may assume that $B = \{2,3,4,5\}$.
		So, by proposition~\ref{cre curve} we obtain 
		that $\Cr(\ell_{1,2}) = \ell_{1,2}$, i.e.
		$\Cr({\mathcal C}_0^{1,2}) = {\mathcal C}_0^{1,2}$.
	\item[(c)] $P_1$ nor $P_2$ is used for the transformation~(\ref{cubic}).
		Then we may assume that $B = \{3,4,5,6\}$.
		So, by proposition~\ref{cre curve} we obtain 
		that $\Cr(\ell_{1,2}) = \ell_3( 3 ; 1^6)$, i.e.
		$\Cr({\mathcal C}_0^{1,2}) = {\mathcal C}_1^{7,8}$.
\end{itemize}
Since we can do this for any $i,j$ 
we obtain that all ${\mathcal C}_1^{i,j}$ are classes of $(-1)$-curves.
Similarly, one can see that, for $a$ odd,
\begin{equation}
	\Cr({\mathcal C}_a^{1,2}) =
		\begin{cases}
			{\mathcal C}_{a+1}^{3,4}  &  \mbox{if } B = \{1,2,3,4\} \\
			{\mathcal C}_a^{1,2}      &  \mbox{if } B = \{2,3,4,5\} \\
			{\mathcal C}_{a-1}^{7,8}  &  \mbox{if } B = \{3,4,5,6\};
		\end{cases}
\end{equation}
and, for $a$ even (and $a >0$),
\begin{equation}
	\Cr({\mathcal C}_a^{1,2}) =
		\begin{cases}
			{\mathcal C}_{a-1}^{3,4}  &  \mbox{if } B = \{1,2,3,4\} \\
			{\mathcal C}_a^{1,2}      &  \mbox{if } B = \{2,3,4,5\} \\
			{\mathcal C}_{a+1}^{7,8}  &  \mbox{if } B = \{3,4,5,6\};
		\end{cases}
\end{equation}
So, we can obtain all classes of type ${\mathcal C}_a^{i,j}$,
and no others.
\end{proof}

\begin{rem}
Lemma~\ref{-1 = C_a^{b,c}} implies that ${\mathcal C}_a^{b,c}$ 
contains precisely one (irreducible) curve, which we denote by $C_a^{b,c}$.
If $a$ is even, $C_a^{b,c}$ is the strict transform of 
a curve of degree $2a+1$ with multiplicity 
$\frac{a}{2}$ at $P_i$ for $i \notin \{b,c\}$ 
and multiplicity $\frac{a}{2}+1$ at $P_b$ and $P_c$.
If $a$ is odd, $C_a^{b,c}$ is the strict transform of 
a curve of degree $2a+1$ with multiplicity 
$\frac{a+1}{2}$ at $P_i$ for $i \notin \{b,c\}$                                        
and multiplicity $\frac{a-1}{2}$ at $P_b$ and $P_c$.   
\end{rem}

\section{Blowings-up of Hirzebruch surfaces along general points}\label{hirz}

Let $\f_n$ be a Hirzebruch surface with $n >0$, 
then $\pic \f_n = \langle f , h_n \rangle = \langle f , c_n \rangle$ with 
$f^2 = 0$, $h_0 . f =1$, $h_n^2 = n$, $c_n = h_n - nf$ and $c_n^2 = -n$
(see e.g.~\cite[Proposition~IV.1]{AB}.

Now let 
$\f^j_n$ be the blowing-up of $\f_n$ along $j$ general points.
By abuse of notation, let $f$, $h_n$ and $c_n$
also denote the pullbacks of these curves on $\f^j_n$,
then $\pic \f^j_n = \langle f , c_n , e_1 , \ldots , e_j \rangle $,
where $e_1, \ldots, e_j$ are the exceptional divisors.

\begin{lemma}\label{F_n met 1 -> F_n-1}
The surface $\f^1_n$, with $n>1$, can also be seen as the blowing-up of
an $\f_{n-1}$ along a general point of $c_{n-1}$.
In particular 
$\pic \f^1_n = \langle f , h_{n-1} , e'_1 \rangle$
with $h_{n-1} = h_n - e_1$
and $e'_1 = f - e_1$ the exceptional divisor corresponding to the blown-up point on $c_{n-1}$.
Moreover 
$ \alpha f + \beta h_n - m e_1  = (\alpha + \beta -m) f + \beta h_{n-1} - (\beta -m) e'_1$.\nvs
\end{lemma}

\begin{proof}
The last equality follows immediately, using $h_{n-1} = h_n - e_1$
and $e'_1 = f - e_1$, so $\pic \f^1_n = \langle f , h_{n-1} , e'_1 \rangle$
is true.
Note that ${e'_1}^2 = -1$ and $h_{n-1}^2 = n-1$.
Consider $c_{n-1} = h_{n-1} - (n-1)f = c_n + e'_1$,
then $c_{n-1}^2 = -(n-1)$.

Now, let $b : \f^1_n \map V$ denote the map obtained by blowing down $e'_1$, 
then $\pic (V) = \langle f , h_{n-1} \rangle = \langle f , c_{n-1} \rangle$,
and $b(e'_1) = Q_1$ is a (general) point on $c_{n-1}$
which is an irreducible curve on $V$ of negative self-intersection. 
So $V = \f_{n-1}$ and $\f^1_n$ is the blowing-up of $\f_{n-1}$ along the point
$Q_1 \in c_{n-1}$.
\end{proof}

\begin{cor}\label{F_n -> P^2}
The surface $\f^{n-1}_n$ can also be seen as the blowing-up of an $\f_1$ 
along $n-1$ general points of $c_1$.
In particular
$\pic \f^{n-1}_n = \langle f , h_{1} , e'_1, \ldots, e'_{n-1} \rangle$
with $h_1 = h_n - e_1 - \cdots - e_{n-1}$ and 
$e'_i = f - e_i$ for all $i =1 ,\ldots,n-1$.
Moreover
$ \alpha f + \beta h_n - m_1 e_1 - \cdots - m_{n-1} e_{n-1}  
	= (\alpha + (n-1) \beta - m_1 - \cdots - m_{n-1}) f + \beta h_{n-1} - (\beta -m_1) e'_1
		- \cdots -(\beta -m_{n-1}) e'_{n-1}$.\nvs
\end{cor}

\begin{proof}
This follows immediately by applying lemma~\ref{F_n met 1 -> F_n-1} $n-1$ times.
\end{proof}

\section{Base locus of linear systems on $X$}\label{results}

A point $P$ of $X$ is called a base point of a linear system
$\cl = \cl_3(d ; m_1,\ldots, m_r)$
if $P \in D$ for all $D \in \cl$.

A divisor $F$ on $X$ is called a fixed component of $\cl$
if $F \subset D$ for all $D \in \cl$.

The base locus of $\cl$, which we denote by $\bs (\cl)$, is defined as
the scheme-theoretical union of all base points.

\begin{exa}
$\bs (\cl_3(2;2^3)) = 2H$, with $H$ the unique element of $\cl_3(1;1^3)$.                 
\end{exa}

Since an empty system obviously has no base locus, 
we only consider non-empty linear systems on $X$
(the results from~\cite{DL} can be used to 
determine wether or not a system is empty).

The main results of this paper are the following

\begin{theor}\label{bs stand}
Let $\cl = \cl_3(d ; m_1,\ldots, m_r) = \cs + \sum_{i=4}^r c_i \cs_i$ be 
(non-empty and) in standard form on $X$,
then the following holds \nvs
\begin{enumerate}
    \item if $d \geq m_1 + m_2$ and $\cl \notin \{\cl_3(2m;m^8),\cl_3(2m;m^7,m-1)\}$
        then $\cl$ is base point free;
    \item if $\cl = \cl_3(2m;m^8)$ ($m \geq 1$) then $\Bs(\cl) = m D_{Q_8}$;
    \item if $\cl = \cl_3(2m;m^7,m-1)$ ($m \geq 1$) then $\Bs(\cl) = mP$
        where $P$ is the unique base point of $\cl|_{Q_8}$
        (which is a point on $D_{Q_8}$);
    \item if $d < m_1 + m_2$,
        then $\Bs(\cl) = \sum_{t_{i,j} > 0} t_{i,j} \ell_{i,j}$,
        with $t_{i,j} = m_i + m_j - d$ ($i \neq j$) and
        $\ell_{i,j}$ the strict transform of the line through $P_i$ and $P_j$.
\end{enumerate}
\end{theor}

\begin{rem}\label{bs stand rem}
Theorem~\ref{bs stand} implies in particular that a class in standard form does not
have fixed components.
\end{rem}

\begin{theor}\label{fixed comp}
Consider the (non-empty) linear system $\cl = \cl_3(d ; m_1,\ldots, m_r)$ on $X$,
then one can obtain the fixed components of $\cl$ as follows \nvs
\begin{itemize}
    \item[(1)]
        Renumber the multiplicities such that $m_1 \geq m_2 \geq \cdots \geq m_r$.
    \item[(2)]
        If $2d < m_1 + m_2 + m_3 + m_4$ then apply the cubic Cremona transformation to
        these 4 multiplicities and goto \begin{rm}(1)\end{rm};
        otherwise goto \begin{rm}(3)\end{rm}.
    \item[(3)]
        If $m_i <0$ then $-m_iE'_i$ is a fixed component,
        and you can apply the cubic Cremona transforms in the opposite direction
        to obtain the class $F_i$ that corresponds to $E'_i$ in the original situation;
        $-m_i F_i$ then belongs to the fixed components of $\cl$.
        Moreover, in this way you obtain all fixed components of $\cl$.
\end{itemize}
\end{theor}

\begin{theor}\label{bs non stand} 
Let $\cl = \cl_3(d ; m_1 , \ldots ,m_r)$ be 
a (non-empty) linear system on $X$ with $m_1 \geq \cdots \geq m_r$.
Assume that $\cl$ has no fixed components and that $\cl$ is not in standard form.
Define $t_a^{b,c} := - \cl . {\mathcal C}_a^{b,c}$,
then the following holds
\begin{itemize}
	\item[(1)] If $4d - \sum_{i=1}^r m_i \neq 1$ then
		$$  \Bs (\cl) = \sum_{t_a^{b,c} > 0} t_a^{b,c} C_a^{b,c}. $$                                                   
	\item[(2)] If $4d - \sum_{i=1}^r m_i (= \cl . D_{Q_8}) = 1$ then $\cl$ can be 
		transformed, by a finite number 
		of Cremona transformations, into $\cl_3(2m ; m^7,m-1)$ for some $m >0$,
		and
		$$  \Bs (\cl) = \sum_{t_a^{b,c} > 0} t_a^{b,c} C_a^{b,c} + mP, $$
		with $P$ the unique base point of $\cl$ on $D_{Q_8}$.
\end{itemize}
\end{theor}

\begin{rem}
Using theorems~\ref{bs stand}, \ref{fixed comp} and \ref{bs non stand},
we can completely determine the base locus of any linear system
$\cl$ on $X$.
\end{rem}

\begin{exa}~\\
Consider the linear system $\cl = \cl_3(15; 13,10,9,7,6,3^2,2)$ on $X$.
\\
First,
we apply the algorithm of theorem~\ref{fixed comp} to determine
the fixed components of $\cl$. 
We use the
following diagram (where Step 1 consists of marking the four
biggest multiplicities):
$$
\begin{array}{c|cccccccc}
 15  &  \framebox{13}  &  \framebox{10}  &  \framebox{9}  &  \framebox{7}  &  6  &  3  &  3  &  2  \\
 6   &  \framebox{4} &  1  &  0  &  -2  &  \framebox{6}  &  \framebox{3}  &  \framebox{3}  & 2  \\
 2   &  \framebox{0} &  \framebox{1}  &  0  &  -2  &  \framebox{2}  &  -1  &  -1  & \framebox{2}  \\
 1   &  -1 &  0  &  0  &  -2  &  1  &  -1  &  -1  & 1  \\
\end{array}
$$
So, after applying the cubic Cremona transform 3 times, we obtain
that $E'_1 + 2 E'_4 + E'_6 + E'_7$ is the fixed part. In order to
go back to the original situation, we now apply the three cubic
Cremona transforms in opposite order. For instance, for $E'_1$ we
obtain
$$
\begin{array}{c|cccccccc}
 0   &  \framebox{-1} &  \framebox{0}  &  0  &  0  &  \framebox{0}  &  0  &  0  & \framebox{0}  \\
 1   &  \framebox{0} &  1  &  0  &  0  &  \framebox{1}  &  \framebox{0}  &  \framebox{0}  & 1  \\
 2   &  \framebox{1}  &  \framebox{1}  &  \framebox{0}  &  \framebox{0}  &  2  &  1  &  1  &  1  \\
 4   &  3 &  3  &  2  &  2  &  2  &  1  &  1  & 1  \\
\end{array} 
$$
Proceeding in the same way for the other $E'_i$, 
we obtain that the fixed components of $\cl$ are
$F = F_1 + 2 F_2 + F_3 + F_4$,
with $F_1 \in \cl_3(4; 3^2 , 2^3 , 1^3)$, 
$F_2 \in \cl_3(1; 1^3)$,
$F_3 \in \cl_3(2; 2, 1^4, 0, 1, 0)$
and $F_4 \in \cl_3(2; 2, 1^5)$.
\\
Now consider $\cl':= \cl - F$, 
then $\cl' = \cl_3(5 ; 4, 3^3, 2, 1^3)$
is a system without fixed components and not in standard form,
so we can apply theorem~\ref{bs non stand} to obtain that
$$
\bs (\cl') = \! \sum_{2 \leq i \leq 4} \! 2 \, C_0^{1,i} + \, C_0^{1,5} + 
	\!\!\!\! \sum_{2 \leq i < j \leq 4} \!\!\!\!\! C_0^{i,j}
    	+ \!\!\!\! \sum_{6 \leq i < j \leq 8} \!\!\!\!\! C_1^{i,j}
$$ 
and
$$
\bs (\cl) =  F_1 + 2 F_2 + F_3 + F_4 + \bs(\cl').
$$
\end{exa}

\section{Proof of theorem~\ref{bs stand}}

Without loss of generality, we may assume that $m_r  >0$.

\begin{rm}(1)\end{rm}
In case $r >4$, we consider the following exact sequence
\begin{equation}\label{ex1}
\xymatrix@1{
	0  \ar[r]  &  \cl - \cs_r  \ar[r]  &  \cl  \ar[r]  &  \cl \otimes \co_{Q_r}  \ar[r]  & 0.
} \end{equation}
We then have that
$\cl \otimes \co_{Q_r}
 = \cl_{Q_r}( d,d; m_1,\ldots,m_r) = \cl_2(2d-m_1 ; (d-m_1)^2 , m_2 ,\ldots,m_r)$
Also, using~\cite[Theorem~5.3]{DL}, we know that
$h^1(\cl - \cs_r) = h^1(\cl) = 0$,
so  
$\cl|_{Q_r} = \cl \otimes \co_{Q_r}$,
i.e. 
$$ \cl|_{Q_r} = \cl_2(2d-m_1 ; (d-m_1)^2 , m_2 ,\ldots,m_r).$$
Since $d \geq m_1 + m_2$, we see that $d - m_1 \geq m_2 (\geq m_3 \geq \cdots \geq m_r)$.
On the other hand $2d - m_1 \geq 2(d-m_1)+m_2$
and 
$\cl|_{Q_r} . K_{Q_r} = -4d + m_1 + \cdots + m_r < -1$
(the inequality is true because $\cl \notin \{\cl_3(2m;m^8),\cl_3(2m;m^7,m-1)\}$).
This means that we can apply \cite[Theorem 3.1 and Corollary 3.4]{BH1}
to conclude that $\cl|_{Q_r}$ is base point free
or thus that $\cl$ has no base points on $Q_r$.

Proceed using the exact sequence~(\ref{ex1}), substituting $\cl$ by $\cl - \cs_r$,
then by $\cl - 2 \cs_r$ and so on,
untill the residue class becomes $\cl - c_r \cs_r$.

Now, let $b$ be $\{ \max \{ i < r : c_i > 0 \}$,
and, if $b \geq 4$, again use the same arguments, now using $Q_b$ in stead of $Q_r$.

Continuing in this way, you can reduce proving the base point freeness of $\cl$ 
to proving the base point freeness of $\cs$.

In order to see that $\cs$ is base boint free,
let $H$ be the unique element of $\cl_3(1 ; 1^3)$ 
and consider the following exact sequence
\begin{equation}\label{ex2}
\xymatrix{
 0  \ar[r]  &  \cl -  \cl_3(1 ; 1^3) \ar[r]  &  \cl  \ar[r]  &  \cl \otimes \co_{H}  \ar[r]  & 0. }
\end{equation}
Using~\cite[Theorem~5.3]{DL}, we know that
$h^1(\cl - \cl_3(1 ; 1^3)) = h^1(\cl) = 0$,
so  
$\cl|_H = \cl \otimes \co_{H}$.
But $\cl \otimes \co_{H} = \cl_2(d - 2m_4 ; m_1-m_4 , m_2-m_4, m_3-m_4)$,
which is base point free (since $d \geq m_1 + m_2$).

Again, we can use this procedure, to see that
$\cs$ is base point free if $\cl_3(d - m_3 - m_4 ; m_1-m_3 , m_2-m_3)$.

Then proceeding in the same way, but use a general $H' \in \cl_3(1 ; 1^2)$
untill the residu class is  $\cl_3(d - m_2 - m_4 ; m_1 - m_2)$;
and after this, using a general $H'' \in \cl_3(1 ; 1)$
untill the residu class is $\cl_3(d - m_1 - m_4)$.

So we actually only need to prove that $\cl_3(d - m_1 - m_4)$ is base point free,
but this is obviously true since $d - m_1 - m_4 \geq 0$.
\vs

\begin{rm}(2)\end{rm}
We use induction on $m$ to prove that $\Bs(\cl) = m D_{Q_8}$.

In case $m=1$, $Q_8 \in \cl = \cs_8$ and we can consider the
following exact sequence
$$
\xymatrix@1{
	0 \ar[r]  &  \co_X  \ar[r]  &  \cs_8 \ar[r]  &  \cs_8 \otimes \co_{Q_8} \ar[r]  &  0. }
$$
Since $h^1(\cs_8)=h^1(\co_X)=0$, we have that 
$\cs_8 |_{Q_8} = \cs_8 \otimes \co_{Q_8}$.
So $\cs_8 |_{Q_8} = \cl_2(3 ; 1^9) = -K_{Q_8}$,
and $D_8$, the unique element of $-K_{Q_8}$,
is the fixed locus of $\cs_8 |_{Q_8}$ and thus also of $\cl$.

Now assume that $m >1$ and that the statement is true for all $m' \leq m-1$.
Consider the exact sequence
$$
\xymatrix@1{
	0  \ar[r]  &  \cl - \cs_8  \ar[r]  &  \cl  \ar[r]  &  \cl \otimes \co_{Q_8}  \ar[r]  & 0.}
$$
Using $h^1(\cl - \cs_8) = h^1(\cl) = 0$ (which follows from~\cite[Theorem~5.3]{DL}),
we obtain that  
$\cl|_{Q_8} = \cl \otimes \co_{Q_8} = -mK_{Q_8}$,
whose only element is $mD_{Q_8}$.
So $\bs (\cl) = D_{Q_8} + \bs(\cl - \cs_8)$,
and, $\cl - \cs_8 = \cl_3(2(m-1); (m-1)^8)$,
so by induction we obtain that $\bs (\cl) = mD_{Q_8}$.
\vs

\begin{rm}(3)\end{rm}
The same procedure as in (2) can be used.
The only difference being that 
$\cl|_{Q_8} = \cl_2(3m ; m^8 , m-1)$,
which has exactly one base point $P$ on $D_{Q_8}$ (see~\cite[Corollary~3.4]{BH1}).
\vs

\begin{rm}(4)\end{rm}
Since obviously $\sum_{t_{i,j} > 0} t_{i,j} \ell_{i,j} \subset \Bs(\cl) $,
it is sufficient to show that $\Bs(\cl) \subset \sum_{t_{i,j} > 0} t_{i,j} \ell_{i,j}$.
To do this, we have to distinguish between $\cs \neq \emptyset$ and $\cs = \emptyset$.

$\bullet$ 
In case $\cs \neq \emptyset$,
$d \geq m_1 + m_4$, and thus also $t_{i,j} \leq 0$ for all $i \geq 1$ and $j \geq 4$.
Now consider the exact sequence~(\ref{ex1}).
Using~\cite[Theorem~5.3]{DL}, we see that 
$h^1(\cl - \cs_r) = h^1(\cl) (\neq 0)$,
and, because of \cite[Lemma~5.2]{DL}, $h^1(\cl \otimes \co_{Q_r}) = 0$,
so 
$\cl|_{Q_r} = \cl \otimes \co_{Q_r} = \cl_2(2d-m_1 ; (d-m_1)^2 , m_2 ,\ldots,m_r)$.
On the other hand, 
$\cl|_{Q_r} . K_{Q_r} = -4d + m_1 + \cdots + m_r < -1$
(the inequality is true because $\cl \notin \{\cl_3(2m;m^8),\cl_3(2m;m^7,m-1)\}$)
and $\cl|_{Q_r}$ is standard (see proof of \cite[Lemma~5.2]{DL}).
This means that we can apply \cite[Theorem 3.1 and Corollary 3.4]{BH1}
to conclude that $\cl|_{Q_r}$ is base point free
or thus that $\cl$ has no base points on $Q_r$.

Continuing this procedure as in (1),
we obtain that $\bs (\cl) \subset \bs (\cs)$.

Now consider the exact sequence~(\ref{ex2}),
then, using \cite[Theorem~5.3]{DL}, 
we obtain that 
$h^1(\cs) = \sum_{t_{i,j} \geq 2} \binom{t_{i,j}+1}{3}$
and
$h^1(\cs -  \cl_3(1 ; 1^3)) = \sum_{t_{i,j} \geq 2} \binom{t_{i,j}}{3}$.
So $h^1(\cs) - h^1(\cs -  \cl_3(1 ; 1^3)) = h^1(\cs \otimes \co_{H})$,
which implies that 
$\cs|_{H} = \cs \otimes \co_{H}$.
Since $\bs (\cs \otimes \co_H) = \sum_{t_{i,j} \geq 1} t_{i,j} \ell_{i,j}$,
we see that $\bs (\cs) = \sum_{t_{i,j} \geq 1} \ell_{i,j} + \bs(\cs -  \cl_3(1 ; 1^3))$.

Again continuing this procedure as in (1),
we finally obtain that $\bs (\cl) = \bs (\cs) = \sum_{t_{i,j} \geq 1} t_{i,j} \ell_{i,j}$.

$\bullet$ 
In case $\cs = \emptyset$, 
$d < m_1 + m_4$, i.e. $t_{1,2} \geq t_{1,3} \geq t_{1,4} > 0$
(and thus also $r \geq 4$). 
Moreover $2d \geq m_1 + m_2 + m_3 + m_4$, so
$d > m_2 + m_3$, and thus
$t_{i,j} \leq 0$ for all $2 \leq i < j$.

Let $W_r$ be a general element of $\cl_3(2 ; 2 , 1^{x})$
with $x = \min \{ r-1 , 5 \}$,
i.e. $W_r$ corresponds in $\p^3$ with an irreducible cone with vertex $P_1$
and through the points $P_2, \ldots, P_{x+1}$.
Then (in $X$) $W_r$ is the blowing-up of a Hirzebruch surface $\f_2$ along $x$ general points,
and $\pic (W_r) = \langle f , h_2 , e_2 , \ldots , e_{x+1} \rangle 
	= \langle f , c_2 , e_2 , \ldots , e_{x+1} \rangle$,
with $c_2 = h_2 - 2 f$,
$\cl_3(1) |_{W_r} = h_2$, $E_1|_{W_r} = c_2$
and $E_i|_{W_r} = e_i$ for all $i = 2,\ldots,x+1$.

Now consider the following exact sequence
$$
\xymatrix@1{
  0  \ar[r]  &  \cl - \cl_3(2 ; 2 , 1^{x})  \ar[r]  &  \cl  \ar[r]  &  
		\cl \otimes \co_{W_r}  \ar[r]  & 0.}
$$
Because of \cite[Theorem~5.3]{DL}, 
we know that 
$h^1(\cl) = \sum_{t_{1,j} > 0} \binom{t_{1,j}+1}{3}$.

\begin{claim}\label{on W_r}
$$ h^1(\cl \otimes \co_{W_r}) = \sum_{\substack{t_{1,j} > 0 \\ j \leq x+1}} \binom{t_{1,j}}{2} \mbox{ and }
	\bs (\cl \otimes \co_{W_r}) = \sum_{\substack{t_{1,j} > 0 \\ j \leq x+1}} t_{1,j} \ell_{1,j}.$$
\end{claim}

\begin{claim}\label{L - W_r}
The linear system $ \cl - \cl_3(2 ; 2 , 1^{x}) $ is in standard form unless 
$\cl = \cl_3(m+m'+t;m'+2t,m',m^6)$ for some $m' \geq m \geq t > 0$. 
Moreover 
$$ h^1(\cl - \cl_3(2 ; 2 , 1^{x}) ) 
	= \sum_{\substack{t_{1,j} > 0 \\ j \leq x+1}} \binom{t_{1,j}}{3} 
		+ \sum_{\substack{t_{1,j} > 0 \\ x+1 < j \leq r}} \! \binom{t_{1,j} +1}{3}. $$
\end{claim}

Using these two claims, we obtain 
that $\cl|_{W_r} = \cl \otimes \co_{W_r}$
and $\bs (\cl|_{W_r}) = \sum_{\substack{t_{1,j} > 0 \\ j \leq x}} t_{1,j} \ell_{1,j}$.
So
$$ \bs (\cl) \subset \bs (\cl - \cl_3(2 ; 2 , 1^{x})) 
	+ \sum_{\substack{t_{1,j} > 0 \\ j \leq x+1}} \ell_{1,j}. $$

If $r >6$, let $H$ be a general element of $\cl_3(1;1,0^5,1^{r-6})$,
denote $\cl - \cl_3(2 ; 2 , 1^{5})$ by $\bar{\cl}$
and consider the following exact sequence
$$
\xymatrix@1{
  0  \ar[r]  &  \cl - \cl_3(3 ; 3 , 1^{r-1})  \ar[r]  &  \bar{\cl} \ar[r]  &  
		\bar{\cl} \otimes \co_{H}  \ar[r]  & 0.}
$$

\begin{claim}\label{on H}
$$ h^1(\bar{\cl} \otimes \co_{H}) 
	= \sum_{\substack{t_{1,j} > 0 \\ 6 < j \leq r}} \binom{t_{1,j}}{2} 
\mbox{ and }
	\bs (\bar{\cl} \otimes \co_{H}) 
		= \sum_{\substack{t_{1,j} > 0 \\ 6 < j \leq r}} t_{1,j} \ell_{1,j}.$$
\end{claim}

\begin{claim}\label{L - W_r - H}
The linear system $ \cl - \cl_3(3 ; 3 , 1^{r-1}) $ is in standard form 
and
$$ h^1(\cl - \cl_3(3 ; 3 , 1^{r-1}) ) 
	= \sum_{t_{1,j} > 0} \binom{t_{1,j}}{3}. $$
\end{claim}

Using claims~\ref{L - W_r} and \ref{L - W_r - H}, we obtain 
that $\bar{\cl}|_{H} = \bar{\cl} \otimes \co_{H}$
and 
$\bs (\bar{\cl}|_{H}) 
	= \sum_{\substack{t_{1,j} > 0 \\ 6 < j \leq r}} t_{1,j} \ell_{1,j}$.

Now define 
$$\cl' = \cl_3(d' ; m'_1 , \ldots, m'_r) := 
	\begin{cases} 
		\cl - \cl_3(2 ; 2 , 1^{x})   &  \mbox{ if } r \leq 6  \\
		\cl - \cl_3(3 ; 3 , 1^{r-1}) &  \mbox{ if } r > 6.
	\end{cases}
$$
and $t'_{i,j} := m'_i + m'_j - d'$,
then
$d' = d-3, m'_1 = m_1 - 3, m'_i = m_i -1 \,\, \forall \, 2 \leq i \leq r,
	\allowbreak t'_{1,j} = t_{1,j} -1 \,\, \forall \, 2 \leq j \leq r$ and 
	$ t'_{i,j} = t_{i,j} +1 \,\, \forall \, 2 \leq i < j \leq r$.
So in particular,
for all $2 \leq i < j \leq r$,
$t'_{i,j} \leq t'_{2,3} = t_{2,3} +1 \leq -t_{1,4} +1 \leq 0$.

If $t'_{1,4} > 0$ (i.e. $t_{1,4} \geq 2$),
then, since $\cl'$ is in standardform (see claims~\ref{L - W_r} and \ref{L - W_r - H}),
we can start our procedure again,
and we can do this
until $t'_{1,4}=0$ for some
$\cl'$
So, in any case, we obtain that
$$ \bs (\cl) \subset \bs (\cl') + \sum_{t_{1,j} > 0} \alpha_{j} \ell_{1,j}, $$
with 
$$
\alpha_{j} =  \begin{cases}
			t_{1,j} & \mbox{ if } m'_j = 0 \\
			\min \{ m_{j} - m'_j, t_{1,j} \}  &  \mbox{ if } m'_j > 0.
		\end{cases}
$$

Because of claims~~\ref{L - W_r} and \ref{L - W_r - H}, 
we also know that $\cl'$ is in standard form,
which means that we are in one of the previously treated cases of our theorem
(since $t'_{1,4}=0$).

If we are in case (1) or in case (4) with $\cs \neq \emptyset$,
then we immediately obtain
$$ \bs (\cl) \subset \sum_{t_{1,j} > 0} t_{1,j} \ell_{1,j}. $$
In case (2), we obtain that
$\cl' = \cl_3( 2m ; m^8)$, for some $m \geq 1$,
$\bs (\cl') = mD_{Q_8}$ and
$\cl = \cl_3( 2m ; m^8) + y \cl_3(3 ; 3 , 1^{7})$ ($y>0$).
So $t_{1,i} = y$ for all $2 \leq i \leq 8$
and
$ \bs (\cl) \subset mD_{Q_8} + \sum^8_{j=1} y \ell_{1,j}$
and, as $D_{Q_8} \subset Q_{8}$,
it is sufficient to prove 
that $\cl$ is base point free on $Q_{8}$.

Consider the exact sequence~(\ref{ex1}).
Then, because of \cite[Theorem~5.3]{DL},
we know that
$h^1(\cl - \cs_8) = h^1(\cl) = 8 \binom{y+1}{3}$.
On the other hand, $h^1(\cl \otimes \co_{Q_8}) = 0$ 
(see \cite[Lemma~5.2]{DL}),
so
$\cl|_{Q_8} = \cl \otimes \co_{Q_8}
	= \cl_2(3m + 3y ; m^2 , (m+y)^7)$,
which is base point free,
since it is in standard form and $\cl|_{Q_8} . K_{Q_8} = -2y \leq -2$ 
(see~\cite[Corollary~3.4]{BH1}).

In case (3), we obtain that
$\cl' = \cl_3( 2m ; m^7, m-1)$, for some $m \geq 1$,
$\bs (\cl') = mP$ and
$\cl = \cl_3( 2m ; m^7, m-1) + y \cl_3(3 ; 3 , 1^{7})$ ($y>0$)
or $\cl = \cl_3( 2 ; 1^7) + y' \cl_3(3; 3, 1^6) + y \cl_3(3 ; 3 , 1^{7})$ ($y,y' \geq 0$ and $y + y'>0$).
Proceeding as above, we can prove that
$\cl|_{Q_8}$ is base point free, and thus obtain that 
$$ \bs (\cl) \subset \sum_{t_{1,j} > 0} t_{1,j} \ell_{1,j}. \bnvs$$
~\vrk

\begin{proof}[Proof of Claim~\ref{on W_r}]
~
\svs \\
We know that 
$\cl \otimes \co_{W_r} = | d h_2 - m_1 c_2 - m_2 e_2 - \cdots -m_{x+1} e_{x+1}|$,
and, because of corollary~\ref{F_n -> P^2},
we obtain
$\cl \otimes \co_{W_r} = |(2d-m_2) h_1 - (d + m_1-m_2) c_1 
	+ t_{1,2} e'_2  - m_3 e_3 - \ldots - m_{x+1} e_{x+1} |$.
So $\cl \otimes \co_{W_r} = t_{1,2} \ell_{1,2} + \cm$
and $\dim (\cl \otimes \co_{W_r}) = \dim (\cm)$,  
with $\cm = \cl_2(2d-m_2 ; d + m_1-m_2 , m_3 ,\ldots, m_{x+1})$.
Using the results of \cite{BH1},
it can be checked that
$$ h^1(\cm) = \sum_{\substack{3 \leq j \leq x+1  \\  t_{1,j} > 0 }} \binom{t_{1,j}}{2}
\,\,\,\, \mbox{ and } \,\,\,\,
\bs (\cm) = \sum_{\substack{3 \leq j \leq x+1  \\  t_{1,j} > 0 }} t_{1,j} \ell_{1,j}.$$
Which then imply that 
$$ h^1(\cl \otimes \co_{W_r}) = \sum_{\substack{j \leq x+1  \\  t_{1,j} > 0 }} \binom{t_{1,j}}{2}
\,\,\,\, \mbox{ and } \,\,\,\, 
\bs (\cl \otimes \co_{W_r}) = \sum_{\substack{j \leq x+1  \\  t_{1,j} > 0 }} t_{1,j} \ell_{1,j}. \bnvs$$
\end{proof}

\begin{proof}[Proof of Claim~\ref{L - W_r}]
~
\svs \\
If $\cl - \cl_3(2 ; 2,1^{x})$ is in standard form, 
then the equality for $h^1(\cl - \cl_3(2 ; 2,1^{x}))$
follows immediately from \cite[Theorem~5.3]{DL}.

Since $d = m_1 + m_4 - t_{1,4}$ and $2d \geq m_1 + m_2 + m_3 + m_4$,
we obtain that 
$m_1 \geq m_2 +  2t_{1,4} + m_3 - m_4 \geq m_2 + 2$.
Using this, it is easy to see
that $2d-4$ is bigger or equal to the sum of the bigget four multiplicities
unless $\cl = \cl_3(m+m'+t;m'+2t,m',m^6)$ for some $m' \geq m \geq t > 0$.
\end{proof}

\begin{proof}[Proof of Claim~\ref{on H}]
~
\svs \\
The statement followd immediately from 
$\bar{\cl} \otimes \co_{H} = \cl_2(d-2 ; m_1-2 , m_7, m_8)$
(or $  = \cl_2(d-2 ; m_1-2 , m_7)$ if $r =7$).
\end{proof}

\begin{proof}[Proof of Claim~\ref{L - W_r - H}]
~
\svs \\
If $\cl - \cl_3(3 ; 3,1^{r-1})$ is in standard form, 
then the equality for $h^1(\cl - \cl_3(3 ; 3,1^{r-1}))$
follows immediately from \cite[Theorem~5.3]{DL}.

So we only need to show that $\cl - \cl_3(3 ; 3,1^{r-1})$ is in standard form,
i.e. that $m_1 - 2 \geq m_2$.
Using $d = m_1 + m_4 - t_{1,4}$ and $2d \geq m_1 + m_2 + m_3 + m_4$,
we obtain that 
$m_1 \geq m_2 +  2t_{1,4} + m_3 - m_4 \geq m_2 + 2$.
\end{proof}

\section{Proof of theorem~\ref{fixed comp}}

Since the Cremona transformation on $X$
is nothing else then blowing-up the lines of
the tetrahedron (formed by the four points used for the transformation) 
and blowing down the other rulings 
of the quadrics obtained in this way (see remark~\ref{cr on X}),
we can eliminate nor construct a fixed part of dimension 2
when applying such a cubic Cremona transformation.

Since we stop applying the Cremona transformation only 
when we obtain something of type
$\cm + \sum m'_i E'_i$
with $m_i > 0$ and $\cm$ a class in standard form,
i.e. a class without fixed components (see remark~\ref{bs stand rem}),
we obtain precisely all fixed components of $\cl$.

\section{Proof of theorem~\ref{bs non stand}}

Proceeding as in the proof of~\cite[Proposition~4.3]{DL}
it is easy to see that 
$F := \sum_{t_a^{b,c} > 0} t_a^{b,c} C_a^{b,c} \subset \bs (\cl)$.
So, if $4d - m_1 - \cdots -m_r \neq 1$ it is enough to prove that there are no base points outside $F$;
and if $4d - m_1 - \cdots -m_r = 1$ it is enough to prove that $\bs (\cl) - F = mP$.

\begin{lemma}\label{L.K = -1 => (3)}
Let $\cl = \cl_3(d; m_1,\ldots,m_r)$ be a (non-empty) class on $X$
which has no fixed components.
Then $4d - m_1 - \cdots -m_r = 1$ if and only if $\cl$ can be 
transformed, by a finite number of Cremona transformations, 
into $\cl_3(2m ; m^7,m-1)$ for some $m >0$.\nvs
\end{lemma}

\begin{proof}
Since the Cremona transformation fixes $D_{Q_8}$ and 
since $\Cr (\cl) . D_{Q_8} = \cl . D_{Q_8}$,
it is clear that $4d - m_1 - \cdots -m_r = 1$ 
if, after a finite number of Cremona transformations, 
$\cl$ transforms into $\cl_3(2m ; m^7,m-1)$.

Conversely, assume $4d - m_1 - \cdots -m_r = 1$.
Then $\cl$ transforms into a class
$\cm = \cl_3(d'; m'_1 , \ldots, m_8)$ in standard form with $m_8 \geq 0$
and $\cm . D_{Q_8} = 4d' - m'_1 - \cdots -m'_8 = 1$,
which implies that $\cm = \cl_3(2m ; m^7,m-1)$.
\end{proof}

\begin{lemma}\label{stand => a=0}
Let $\cn := \cl_3(d; m_1,\ldots,m_8)$ be a class in standard form on $X$,
then $\cn . {\mathcal C}_a^{b,c} \geq 0$ for all $a>0$ and for all $b,c \in \{1,\ldots,8\}$.\nvs 
\end{lemma}

\begin{proof}
If $a$ is even, then $\cn . {\mathcal C}_a^{b,c} \geq \cn . {\mathcal C}_a^{1,2}$
for all $b,c \in \{1,\ldots,8\}$,
and 
$\cn . {\mathcal C}_a^{1,2} = - \frac{a}{2} ( t_0^{1,2} + t_0^{3,4} + t_0^{5,6}) 
	- (\frac{a}{2} -1)t_0^{7,8} -( t_0^{1,2} + t_0^{7,8} )$.
Since $\cn$ is in standard form,
$0 \geq t_0^{3,4} \geq t_0^{5,6} \geq t_0^{7,8}$
and $0 \geq t_0^{1,2} + t_0^{3,4} \geq t_0^{1,2} + t_0^{7,8}$,
so $\cn . {\mathcal C}_a^{1,2} \geq 0$.

If $a$ is odd, then $\cn . {\mathcal C}_a^{b,c} \geq \cn . {\mathcal C}_a^{7,8}$
for all $b,c \in \{1,\ldots,8\}$,
and 
$\cn . {\mathcal C}_a^{7,8} = - \frac{a-1}{2}( t_0^{1,2} + t_0^{3,4} + t_0^{5,6}+ t_0^{7,8} ) 
	- (t_0^{1,2} + t_0^{3,4} + t_0^{5,6}) \geq 0$.
\end{proof}

To simplify notation, assume we want to apply the Cremona transformation using $P_1, \ldots, P_4$.

Let $Y$ be the blowing-up of $X$ along the $\ell_{i,j}$, $1 \leq i < j \leq 4$,
$p: Y \map X$ the projection map,
let $E_i$, $F_i$, $E_{i,j}$ and $F_{i,j}$ be as in (\ref{pic Y}) and (\ref{F i,j})
and let $h$, $h'$, $e_i$, $f_j$, $e_{i,j}$ and $f_{i,j}$ be as in (\ref{A^2(Y)}) and (\ref{f i,j}).

Let $p' : Y \map X'$ be the map obtained by blowing down the $F_{i,j}$.

Now, analoguously to $C_a^{b,c}$, define $D_a^{b,c}$
in ${\rm{\bf A}}^2(X')$
(e.g. $D_1^{7,8} = | 3 h' - f_1 - \cdots - f_4 - e_5 - e_6 |$).
And define $s_a^{b,c} := - \Cr (\cl) . D_a^{b,c}$.

By abuse of notation, if $a > 0$ or if $a = 0$ and $\{b,c\} \not\subset \{1,2,3,4\}$,
we also denote the pull-back of $C_a^{b,c}$, resp. $D_a^{b,c}$, by $C_a^{b,c}$, resp. $D_a^{b,c}$.

Let $F^*$ denote the pull-back on $Y$ of $F$,
and write $F^*$ as $F^{(1)} + F^{(2)}$,
with  
$$
F^{(2)} = \sum_{\substack{1 \leq b < c \leq 4 \\ t_0^{b,c} > 0}} \!\!\!\! t_0^{b,c} E_{b,c}.
$$

Similarly, let $G^*$ denote the pull-back
of $G = \sum_{s_a^{b,c} > 0} s_a^{b,c} D_a^{b,c}$ on $Y$,
and write $G^*$ as $G^{(1)} + G^{(2)}$,
with  
$$
G^{(2)} = \sum_{\substack{1 \leq b < c \leq 4 \\ s_0^{b,c} > 0}} \!\!\!\! s_0^{b,c} F_{b,c}.
$$

Define  
$ \cm := p^*(\cl) \otimes \co_{Y}(-F^{(2)}) \otimes \ci_{F^{(1)}}$.

\begin{prop}\label{cr(M)}
$$ \cm = {p'}^*(\Cr (\cl)) \otimes \co_{Y}(-G^{(2)}) \otimes \ci_{G^{(1)}}.  \nvs $$
\end{prop}

\begin{proof}
First of all, by abuse of notation, let us write
$\cm$ as $\cm^{(2)} - \cm^{(1)}$,
with 
$$ 
\begin{gathered}
\cm^{(2)} = dH - \!\! \sum_{1\leq i \leq 8} \!\! m_i E_i - 
     \!\!\!\!\!  \sum_{\substack{1 \leq i < j \leq 4 \\ t_0^{i,j} >0}} \!\!\!\!\! 
			t_0^{i,j} E_{i,j} 
\\
\!\!\!\!\!\!\!\!\!\! \mbox{and } \,\,\,
\cm^{(1)} = \!\!\!\!\!\!\!\!\!\!\!  
	\sum_{\substack{a > 0 \mbox{\scriptsize{ or }} a=0 \mbox{\scriptsize{ and }} 4<c \\ t_a^{b,c} >0}} 
	\!\!\!\!\!\!\!\!\!\!\! 	t_a^{b,c} C_a^{b,c}. 
\end{gathered}
$$
Using the formulas~\ref{F i,j}
and the fact that $s_0^{i,j} = d - m_k - m_l = -t_0^{k,l}$ with
$\{i,j,k,l\} = \{1,2,3,4\}$,
it can easily be checked that, if $s = 2d - \sum_{i=1}^4 m_i$,
$$ 
\cm^{(2)} = (d+s) H' - \sum_{1 \leq i \leq 4} (m_i + s) F_i - \sum _{5 \leq i \leq r} m_i E_i
	- \!\!\!\! \sum_{\substack{1 \leq i < j \leq 4 \\ s_0^{i,j} > 0}} \!\!\! s_0^{i,j} F_{i,j}. 
$$
Moreover, using the formulas~\ref{f i,j},
a simple calculation shows that 
$$ 
\cm^{(1)} = \!\!\!\!\!\!\!\!\!\!\!  
	\sum_{\substack{a > 0 \mbox{\scriptsize{ or }} a=0 \mbox{\scriptsize{ and }} 4<c \\ s_a^{b,c} >0}} 
	\!\!\!\!\!\!\!\!\!\!\! 	s_a^{b,c} D_a^{b,c}. 
$$
Combining these two results,
we obtain that
$$\cm = {p'}^*(\Cr (\cl)) \otimes \co_{Y}(-G^{(2)}) \otimes \ci_{G^{(1)}}.$$
\end{proof}

\begin{cor}\label{bs L -> bs Cr L}
$\bs (\cl) - F \neq \emptyset$ if and only if $\bs (\Cr(\cl)) - G \neq \emptyset$.\nvs
\end{cor}

\begin{proof}
Since the Cremona transformation is an involution, 
it is sufficient to prove just one implication.
If $P \in \bs (\cl) - F$ then 
$p^{-1}(P) \subset \bs (\cm) - F^{*} = \bs (\cm) - G^{*}$,
which implies that $p'(p^{-1}(P)) \subset  \bs (\Cr(\cl)) - G$
(and  $p'(p^{-1}(P)) \neq \emptyset$).
\end{proof}

(1) 
If $4d - m_1 - \cdots - m_r \neq 1$,
we apply Cremona until we obtain a class 
$\cl'$ in standard form.
Because of corollary~\ref{bs L -> bs Cr L}, it is enough to prove that
$\bs (\cl') - F' = \emptyset$,
with $F' = \sum_{{t'}_a^{b,c} > 0}{t'}_a^{b,c} {C'}_a^{b,c}$.
But, because of lemma~\ref{stand => a=0}, ${t'}_a^{b,c} \leq 0$ if $a >0$,
i.e. we obtain that $F' = \sum_{{t'}_0^{b,c} > 0}{t'}_0^{b,c} {C'}_0^{b,c}$.
On the other hand, 
$4d - m_1 - \cdots - m_r \neq 1$ implies that 
$\cl'$ is of type (1) or (4) of theorem~\ref{bs stand}
($\cl'$ is not of type (2) since this is a class which is invariant under Cremona,
and it is not of type (3) because of lemma~\ref{L.K = -1 => (3)}).
So it follows from theorem~\ref{bs stand}
that $\bs (\cl') = F'$, and thus
$\bs (\cl') - F' = \emptyset$.

(2)
If $4d - m_1 - \cdots - m_r = 1$,
we apply Cremona until we obtain the class 
$\cl' = \cl_3(2m ; m^7 , m-1)$ (see lemma~\ref{L.K = -1 => (3)}).
Reasoning as before,
$F' = \sum_{{t'}_0^{b,c} > 0}{t'}_0^{b,c} {C'}_0^{b,c}$,
but now ${t'}_0^{b,c}$ is either equal to 0 or -1,
so $F' = \emptyset$.
On the other hand,
because of theorem~\ref{bs stand},
$\bs (\cl') = mP'$,
and, since $P'$ is never on a strict transform of an edge of the tetrahedron used for the Cremona transformation,
proceeding as in the proof of corollary~\ref{bs L -> bs Cr L},
we obtain that, on $X$,
$P'$ corresponds to the base point $P$ of $\cl_3(2; 1^7)$ on $D_{Q_8}$.
So we obtain that
$\bs (\cl) - F = mP$.



\end{document}